\documentclass[twoside,reqno,A4]{amsart}

\usepackage{latexsym}

\newcommand{\é}{\'e}

\newcommand{\qdp}{\hspace*{1.5mm}}%

\newcommand{\xquad}{\quad\qquad\quad}
\newcommand{\qdn}{\hspace*{-1.5mm}}
\newcommand{\qqdn}{\hspace*{-2.5mm}}
\newcommand{\xqdn}{\hspace*{-5.0mm}}

\newcommand{\sst}{\scriptstyle}
\newcommand{\sss}{\scriptscriptstyle}

\newcommand{\+}{&\qqdn}%
\newcommand{\zero}{&\qdn&\xqdn}
\newcommand{\tagg}[2]{\addtocounter{equation}{#1}
		\tag{\theequation#2}}

\newcommand{\Lra}{\Longrightarrow}

\newcommand{\ang}[1]{\langle{#1}\rangle}

\newcommand{\mc}[1]{\mathcal{#1}}
\newcommand{\sbs}[1]{_{\substack{#1}}}

\newcommand{\binm}{\binom}
\newcommand{\sbnm}[2]{\Bigl(\!\ba{c}\!#1\!\\#2\ea\!\Bigr)}

\newcommand{\double}[2]{\genfrac{}{}{0mm}{1}{#1}{#2}}

\newcommand{\be}{\begin{equation}}
\newcommand{\ee}{\end{equation}}
\newcommand{\ba}{\begin{array}}
\newcommand{\ea}{\end{array}}
\newcommand{\bmn}{\begin{eqnarray}}
\newcommand{\emn}{\end{eqnarray}}
\newcommand{\bnm}{\begin{eqnarray*}}
\newcommand{\enm}{\end{eqnarray*}}
\newcommand{\bln}{\begin{subequations}}
\newcommand{\eln}{\end{subequations}}
\newcommand{\blm}{\bln\bmn}
\newcommand{\elm}{\emn\eln}

\newcommand{\eqn}[1]{\begin{equation}#1\end{equation}}
\newcommand{\xalignx}[1]{\begin{align}#1
            \end{align}}     
\newcommand{\xalignz}[1]{\begin{align*}#1
            \end{align*}}    

\newcommand{\mult}[2]{\begin{array}{#1}#2\end{array}}%

\newcommand{\lam}{\lambda}

\newcommand{\del}{\delta}

\newcommand{\sig}{\sigma}
\newcommand{\ome}{\omega}

\newcommand{\Ome}{\Omega}

\newtheorem{thm}{Theorem}

\newtheorem{corl}[thm]{Corollary}

\newtheorem{exam}{Example}

\newtheorem{entry}{Entry}

\newcommand{\bbtm}[4]{\bibitem{kn:#1}{#2,}~\emph{#3,}~{#4.}}	
\newcommand{\cito}[1]{\cite{kn:#1}}	
\newcommand{\citu}[2]{\cite[#2]{kn:#1}}

\begin{document} 

\title{Partial-Fraction Decompositions\\
	 and Harmonic Number Identities}
\author{CHU Wenchang}
\address{\textbf{Prof.\:CHU Wenchang}	\newline
 	Dipartimento\: \:di\: \:Matematica	\newline
     	Universit\a`{a} degli Studi di Lecce\newline
     	Lecce-Arnesano\: \:P.\:O.\:Box\:193	\newline
      	73100 Lecce,$\quad$ITALIA	\newline              	
	tel 39+0832+297409			\newline 
	fax 39+0832+297594			\newline              	
	Email \emph{chu.wenchang@unile.it}}

\date{\today,$\quad$The work carried out during my visit 
to Center for Combinatorics (LPMC), Nankai University (2005)}%


\begin{abstract}
By means of partial fraction method, we investigate 
the decomposition of rational functions.
Several striking identities on harmonic numbers 
and generalized Ap\éry numbers will be established,
including the binomial-harmonic number identity 
associated with Beukers' conjecture on Ap\éry numbers. 
\end{abstract} 

\maketitle
\markboth{CHU Wenchang: Partial-Fraction Decompositions 
                        and Harmonic Number Identities}
         {CHU Wenchang: Partial-Fraction Decompositions
                        and Harmonic Number Identities}
\thispagestyle{empty}


\section{Introduction}

The generalized harmonic numbers are defined 
to be partial sums of the Riemann-Zeta series: 
\eqn{H_0^{\ang{m}}=0
\quad\text{and}\quad 
H_n^{\ang{m}}=\sum_{k=1}^n
\frac{1}{k^m}
\quad\text{for}\quad 
m,\:n=1,2,\cdots.}
When $m=1$, they reduce to the classical 
ones, shortened as $H_n=H_n^{\ang{1}}$.

If the shifted factorial is defined by  
\eqn{(c)_0\equiv1
\quad\text{and}\quad
(c)_n=c(c+1)\cdots(c+n-1)
\quad\text{for}\quad
n=1,2,\cdots}
then we can establish, by means of the standard 
partial-fraction decompositions, the following 
algebraic identities:
\xalignz{
\frac{n!}{(x)_{n+1}}
\:=\:&\tagg{1}{}\label{pf0/1}
\sum^n_{k=0}
\sbnm{n}{k}
\frac{(-1)^k}{x+k}\\
\frac{(n!)^2}{(x)^2_{n+1}}
\:=\:&\tagg{1}{}\label{pf0/2}
\sum^n_{k=0}
\sbnm{n}{k}^2
\Big\{\frac{1}{(x+k)^2}+
\frac{2}{x+k}(H_k-H_{n-k})\Big\}\\
\frac{(n!)^3}{(x)^3_{n+1}}
\:=\:&\tagg{1}{a}\label{pf0/3a}
\sum^n_{k=0}
(-1)^k\sbnm{n}{k}^3
\bigg\{\frac{1}{(x+k)^3}+
\frac{3}{(x+k)^2}(H_k-H_{n-k})\\
&\;+\:\tagg{0}{b}\label{pf0/3b}
\frac{3}{2(x+k)}
\Big[3(H_k-H_{n-k})^2
+\big(H_k^{\ang{2}}+H^{\ang{2}}_{n-k}
\big)\Big]\bigg\}.}
Multiplying (\ref{pf0/3a}-\ref{pf0/3b}) across 
by $x$ and then letting $x\to\infty$, we recover 
one identity among the hardest challenges
claimed in \citu{pade-2}{Eq\:16}, \citu{pade-3}{Eq\:12}
and \citu{pade-1}{Eq\:20}:
\eqn{\sum_{k=0}^n(-1)^k\sbnm{n}{k}^3
\Big\{3(H_k-H_{n-k})^2
+\big(H_k^{\ang{2}}+H^{\ang{2}}_{n-k}\big)\Big\}
\:=\:0.\label{pf0/3c}}

This has best exemplified the power of partial fraction method. 
For more general rational functions, we will investigate their 
partial fraction decompositions in the second section, which 
involve the complete Bell polynomials (or cyclic indicators 
of symmetric groups) on the generalized harmonic numbers. 
Several further examples and miscellaneous formulae will be 
collected in the third and last section. In order to facilitate 
consultation for readers, three short tables of the complete 
Bell polynomials on the generalized harmonic numbers will be 
presented in the appendices. 
    

\section{Partial Fraction Decompositions}

For two natural numbers $n$ and $k$ with $0\le k\le n$, 
define two functions related to harmonic numbers by
\blm
H_\ell(x):=\sum_{\iota=1}^n
\frac{1}{(\iota-x)^\ell}
&\Lra&\label{replace-a}
H_\ell(-k)=H^{\ang{\ell}}_{n+k}-H^{\ang{\ell}}_k\\
\mc{H}_\ell(x)
:=\sum\sbs{\iota=0\\\iota\not=k}^n
\frac{1}{(\iota+x)^\ell}
&\Lra&\label{replace-b}
\mc{H}_\ell(-k)=H_{n-k}^{\ang{\ell}}
+(-1)^{\ell}H_k^{\ang{\ell}}.
\elm
They come respectively from the logarithmic 
derivatives of the binomial coefficients 
\blm
h(x)&=&\frac{(1-x)_n}{n!}
\:=\:\sbnm{n-x}{n}\\
\hbar(x)&=&\frac{n!\times(x+k)}{(x)_{n+1}}
\:=\:\frac{\binm{n}{k}}
{\binm{x+k-1}{k}\binm{x+n}{n-k}}.
\elm
  
Let $\sig(\ell)$ be the set of partitions 
of $\ell$ represented by $\ell$-tuples 
of nonnegative integers $({m_1},{m_2},\cdots,{m_\ell})$ 
such that $\sum_{k=1}^\ell km_k=\ell$. 
Its subset of $\ell$-partitions into $m$ parts 
with $\sum_{k=1}^\ell m_k=m$ is denoted 
by $\sig_m(\ell)$. 
\begin{thm}[Partial fraction decomposition]
Let $\lam,\:\mu$ and $n$ be three natural numbers 
with $\lam+(\lam-\mu)n>0$.
Then there holds the algebraic identity:
\eqn{\frac{(n!)^{\lam-\mu}(1-x)_n^\mu}
       {(x)^\lam_{n+1}}
\:=\:\sum^n_{k=0}(-1)^{k\lam}
\sbnm{n}{k}^\lam\sbnm{n+k}{k}^\mu
\sum_{\ell=0}^{\lam-1}
\frac{\Ome_\ell(\lam,\mu,-k)}
     {\ell!\:(x+k)^{\lam-\ell}} }
with the $\Ome$-coefficients being determined 
by the Bell polynomials \emph{(or the cyclic 
indicators of symmetric groups)}: 
\eqn{\Ome_\ell(\lam,\mu,x)
\:=\:\label{bell-x}
(-1)^\ell\ell!
\sum\sbs{\sig(\ell)}
\prod_{i=1}^\ell
\frac{\Big\{\lam\mc{H}_i(x)-(-1)^i\mu{H}_i(x)\Big\}^{m_i}}
     {m_i!\:i^{m_i}}}
where the multiple sum runs over $\sig(\ell)$, the set 
of $\ell$-partitions represented by $\ell$-tuples 
of nonnegative integers $({m_1},{m_2},\cdots,{m_\ell})$ 
such that $\sum_{k=1}^\ell km_k=\ell$.
\end{thm}
In particular, the $\Ome$-coefficients read explicitly as
\eqn{\Ome_\ell(\lam,\mu,-k)
\:=\:\label{bell-k}
\ell!
\sum\sbs{\sig(\ell)}
\prod_{i=1}^\ell
\frac{\Big\{\lam\big[H_k^{\ang{i}}+(-1)^iH_{n-k}^{\ang{i}}\big]
      +\mu\big[H_k^{\ang{i}}-H_{n+k}^{\ang{i}}\big]\Big\}^{m_i}}
     {m_i!\:i^{m_i}}.}
\begin{proof} By means of partial fraction decomposition, 
we can formally write 
\[\frac{\hbar^\lam(x)h^\mu(x)}{(x+k)^\lam}
=\frac{(n!)^{\lam-\mu}(1-x)_n^\mu}
      {(x)^\lam_{n+1}}
\:=\:\sum^n_{k=0}
\sum_{\ell=0}^{\lam-1}
\frac{C(k,\ell)}{(x+k)^{\lam-\ell}}\]
where the coefficients $C(k,\ell)$ are to be determined. 
Letting $\mc{D}_x=\frac{d}{dx}$ stand for the derivative 
operator with respect to $x$ and then noting that
\bnm
\hbar(-k)&=&(-1)^{k}\sbnm{n}{k}\\
h(-k)&=&\sbnm{n+k}{k}
\enm
we first demonstrate that 
for $0\le\ell<\lam$ there holds:
\eqn{C(k,\ell)
\:=\:\label{coeff-a}
\hbar^\lam(-k)h^\mu(-k)
\times\frac{\Ome_\ell(\lam,\mu,-k)}{\ell!}}
where the $\Ome$-coefficients are given 
by the following logarithmic derivatives:
\eqn{\Ome_\ell(\lam,\mu,x)
\quad=\quad\label{coeff-b}
\frac{\mc{D}_x^{\ell}
\Big\{\hbar^\lam(x)h^\mu(x)\Big\}}
     {\hbar^\lam(x)h^\mu(x)}.}

For $\ell=0$, we have obviously 
$\Ome_0(\lam,\mu,x)\equiv1$ and that 
\bnm
C(k,0)&=&\lim_{x\to-k}
\hbar^\lam(x)h^\mu(x)
\times\Ome_0(\lam,\mu,x)\\
&=&\hbar^\lam(-k)h^\mu(-k)
\times\Ome_0(\lam,\mu,-k).
\enm
Next for $\ell=1$, we can check 
(\ref{coeff-a}-\ref{coeff-b})
through L'Hospital's rule that 
\bnm
C(k,1)&=&\lim_{x\to-k}
(x+k)^{\lam-1}
\Big\{\frac{\hbar^\lam(x)h^\mu(x)}{(x+k)^{\lam}}
-\frac{C(k,0)}{(x+k)^{\lam}}\Big\}\\
&=&\lim_{x\to-k}
\frac{\hbar^\lam(x)h^\mu(x)-C(k,0)}{x+k}\\
&=&\lim_{x\to-k}
\mc{D}_x\Big\{\hbar^\lam(x)h^\mu(x)\Big\}\\
&=&\hbar^\lam(-k)h^\mu(-k)
\times\Ome_1(\lam,\mu,-k).
\enm

Supposing now the truth of (\ref{coeff-a}-\ref{coeff-b})
for $\ell=0,1,\cdots,m-1$ with $m<\lam$, then we have 
to verify it also for $\ell=m$. Applying again the L'Hospital 
rule for $m$-times, we can determine the coefficient 
\bnm
C(k,m)&=&\lim_{x\to-k}
(x+k)^{\lam-m}
\Big\{\frac{\hbar^\lam(x)h^\mu(x)}{(x+k)^{\lam}}
-\sum_{\ell=0}^{m-1}
\frac{C(k,\ell)}{(x+k)^{\lam-\ell}}\Big\}\\
&=&\lim_{x\to-k}
\frac{1}{(x+k)^{m}}
\Big\{\hbar^\lam(x)h^\mu(x)
-\sum_{\ell=0}^{m-1}C(k,\ell)
\times(x+k)^{\ell}\Big\}\\
&=&\lim_{x\to-k}\hbar^\lam(x)h^\mu(x)
\frac{\mc{D}_x^m\Big\{\hbar^\lam(x)h^\mu(x)\Big\}}
     {m!\:f^{\lam}(x)}\\
&=&\hbar^\lam(-k)h^\mu(-k)
\times\frac{\Ome_m(\lam,\mu,-k)}{m!}.
\enm
Based on the induction principle, 
we have confirmed that the coefficients
in partial fraction decomposition are
determined by (\ref{coeff-a}-\ref{coeff-b}). 

To complete the proof of the theorem,
it remains to show that these coefficients 
can be calculated explicitly through equation 
\eqref{bell-x} and therefore \eqref{bell-k}
(Bell polynomials and/or the cyclic indicators 
of symmetric groups). 

Manipulating the differential operation
\[\frac{\mc{D}_x^{1+\ell}
\Big\{\hbar^\lam(x)h^\mu(x)\Big\}}
       {\hbar^\lam(x)h^\mu(x)}
\:=\:
\frac{\mc{D}_x
\Big\{\hbar^\lam(x)h^\mu(x)\Big\}}
     {\hbar^\lam(x)h^\mu(x)}
\frac{\mc{D}_x^{\ell}
\Big\{\hbar^\lam(x)h^\mu(x)\Big\}}
     {\hbar^\lam(x)h^\mu(x)}\]
we can derive for \eqref{coeff-b} the recurrence relation 
\eqn{\Ome_{1+\ell}(\lam,\mu,x)
\:=\:\label{recurrence}
\Big\{\mc{D}_x-\lam\mc{H}_1(x)-\mu{H}_1(x)\Big\}
\:\Ome_\ell(\lam,\mu,x).}
It is trivial to see that $\Ome_\ell(\lam,\mu,x)$ 
defined by \eqref{bell-x} admits the initial condition
$\Ome_0(\lam,\mu,x)\equiv1$. If we can check that 
$\Ome_\ell(\lam,\mu,x)$ defined by \eqref{bell-x}
satisfies the same recurrence relation \eqref{recurrence}, 
then the validity of \eqref{bell-x} would be confirmed 
for all the natural numbers $\ell$.

Now substituting the RHS of \eqref{bell-x} 
into the RHS of \eqref{recurrence} and 
then noticing the differential relations 
\bnm
\mc{D}_x {H}_j(x)
&=&+j{H}_{j+1}(x)\\
\mc{D}_x \mc{H}_j(x)
&=&-j\mc{H}_{j+1}(x)
\enm  
we get the following expression
\bln\xalignx{
(&-1)^{1+\ell}\ell!\label{part0}
\Bigg\{\Big[\lam\mc{H}_1(x)+\mu{H}_1(x)\Big]
\sum\sbs{\sig(\ell)}\prod_{i=1}^\ell
\frac{\Big\{\lam\mc{H}_i(x)
	-(-1)^i\mu{H}_i(x)\Big\}^{m_i}}
     {m_i!\:i^{m_i}}\\
&\label{parts}+
\sum\sbs{\sig(\ell)}
\prod_{i=1}^\ell
\frac{\Big\{\lam\mc{H}_i(x)-(-1)^i\mu{H}_i(x)\Big\}^{m_i}\qdn}
     {m_i!\:i^{m_i}}\qdn
\sum_{j=1}^\ell jm_j
\tfrac{\lam\mc{H}_{1+j}(x)+(-1)^j\mu{H}_{j+1}(x)}
      {\lam\mc{H}_j(x)-(-1)^j\mu{H}_j(x)}\Bigg\}.
}\eln

In accordance with the combinatorial structure, each $\ell$-partition 
enumerated by $\sig(\ell)$ becomes a $(1+\ell)$-partition 
with a ``$j$"-part being shifted to a ``$1+j$"-part for $0\le j\le\ell$.
Vice versa, every $(1+\ell)$-partition enumerated by $\sig(1+\ell)$
reduces to a $\ell$-partition with a ``$1+j$"-part being replaced 
by a ``$j$"-part for $0\le j\le\ell$. Then the sum over partitions 
should be reformulated accordingly.
 
First, the line \eqref{part0} with an extra part ``$1$" yields 
a new factor $M_1:=m_1+1$. Then if $m_\ell=0$, for each $j$ 
corresponding to the shift from part ``$j$" to part ``$1+j$" 
displayed in line \eqref{parts}, the coefficient $jm_j$ is 
replaced by $(1+j)M_{j+1}$ under two index substitution 
$M_{j}:=m_j-1$ and $M_{1+j}:=m_{1+j}+1$ for $1\le j<\ell$.
Lastly if $m_\ell=1$, the coefficient $\ell m_\ell$ will 
be replaced by $(1+\ell)M_{\ell+1}$ with two summation index 
being substituted by $M_{\ell}:=m_{\ell}-1$ and $M_{1+\ell}:=m_{\ell}$. 
Summing up, we may combine \eqref{part0} with \eqref{parts} 
and obtain the following expression 
\[(-1)^{1+\ell}\ell!
\sum\sbs{\sig(1+\ell)}
\prod_{i=1}^{1+\ell}
\frac{\Big\{\lam\mc{H}_i(x)
	-(-1)^i\mu{H}_i(x)\Big\}^{M_i}\qdn}
     {M_i!\:i^{M_i}}\qdn
\sum_{j=1}^{1+\ell}jM_j.\]
According to \eqref{bell-x}, the last expression 
becomes $\Ome_{1+\ell}(\lam,\mu,x)$, i.e., 
the left member of \eqref{recurrence} 
thanks for the $(1+\ell)$-partition 
$1+\ell=\sum_{j=1}^{1+\ell}jM_j$.
This confirms that the RHS of \eqref{bell-x}
satisfies indeed the recurrence relation 
\eqref{recurrence}. We therefore have established 
equality \eqref{bell-x} and \eqref{bell-k}.
This completes the proof of the theorem.
\end{proof}

In the theorem, multiplying the partial fraction 
decomposition by $x$ and then letting $x\to\infty$, we derive
the following harmonic number identity.
\begin{corl}[Harmonic number identity]
Let $\lam,\:\mu$ and $n$ be three natural numbers 
with $\lam+(\lam-\mu)n>1$.
Then there holds the algebraic identity:
\eqn{\sum^n_{k=0}(-1)^{k\lam}
\sbnm{n}{k}^\lam\sbnm{n+k}{k}^\mu
{\Ome_{\lam-1}(\lam,\mu,-k)}
\:=\:0 }
where the $\Ome$-coefficients are given 
by the Bell polynomials \eqref{bell-x}.
\end{corl}

If defining further two sequences by
\blm
\varpi_{\ell}(\lam,x)
&=&\Ome_\ell(\lam,0,x)
\:=\:\frac{\mc{D}_x^{\ell}
\hbar^\lam(x)}{\hbar^\lam(x)}
\:=\:\ell!\label{bell-xa}
\sum\sbs{\sig(\ell)}
{\:(-1)^\ell\:\lam^m}
\prod_{i=1}^\ell
\frac{\mc{H}_i^{m_i}(x)}
     {m_i!\:i^{m_i}}\\
\ome_{\ell}(\mu,x)
&=&\Ome_\ell(0,\mu,x)
\:=\:\frac{\mc{D}_x^{\ell}h^\mu(x)}{h^\mu(x)}
\:=\:\ell!\label{bell-xb}
\sum\sbs{\sig(\ell)}
{(-1)^m\mu^m}
\prod_{i=1}^\ell
\frac{H_i^{m_i}(x)}
     {m_i!\:i^{m_i}}
\elm
and then applying the Leibniz rule to \eqref{coeff-b},
we find the following convolution formula:
\eqn{\Ome_\ell(\lam,\mu,x)
\:=\:\label{leibniz-x}
\sum^{\ell}_{\iota=0}
\sbnm{\ell}{\iota}
\varpi_{\iota}(\lam,x)
\ome_{\ell-\iota}(\mu,x).}

Putting $x=-k$, we write down the corresponding 
relation as follows:
\eqn{\Ome_\ell(\lam,\mu,-k)
\:=\:\label{leibniz-k}
\sum^{\ell}_{\iota=0}
\sbnm{\ell}{\iota}
\varpi_{\iota}(\lam,-k)
\ome_{\ell-\iota}(\mu,-k)}
where $\varpi$ and $\ome$ are explicitly 
provided by the following formulae:
\blm
\varpi_{\ell}(\lam,-k)
&=&\Ome_\ell(\lam,0,-k)
\:=\:\ell!\label{bell-ka}
\sum\sbs{\sig(\ell)}{\lam^m}
\prod_{i=1}^\ell
\frac{\Big\{H_k^{\ang{i}}
	+(-1)^iH_{n-k}^{\ang{i}}\Big\}^{m_i}}
     {m_i!\:i^{m_i}}\\
\ome_{\ell}(\mu,-k)
&=&\Ome_\ell(0,\mu,-k)
\:=\:\ell!\label{bell-kb}
\sum\sbs{\sig(\ell)}{\mu^m}
\prod_{i=1}^\ell
\frac{\Big\{H_k^{\ang{i}}
	-H_{n+k}^{\ang{i}}\Big\}^{m_i}}
     {m_i!\:i^{m_i}}.
\elm


\section{Examples: Harmonic Number Identities}

By means of the theorem and the corollary, we will display 
several examples of partial fraction decompositions and 
the corresponding harmonic number identities.

\begin{exam}[$\lam=1$] For $\del=0,\:1$, 
there hold partial fraction expansions
\[(n!)^{1-\del}
\frac{(1-x)_n^{\mu}}{(x)_{n+1}}
\:=\:\sum^n_{k=0}
\sbnm{n}{k}\sbnm{n+k}{k}^{\del}
\frac{(-1)^k}{x+k}\]
and the two corresponding harmonic number identities:
\[(-1)^n\del\:=\:\sum^n_{k=0}(-1)^k
\sbnm{n}{k}\sbnm{n+k}{k}^{\del}.\]
\end{exam}

When $\mu=0$, the corresponding partial fraction
expansion reads as the formula displayed in 
\eqref{pf0/1}. For an alternative derivation, 
refer to the recent paper \cito{koepf}.
When $\mu=1$, the last binomial identity is a special 
case of the Chu-Vandermonde convolution formula
on binomial coefficients:
\[\sbnm{u+v}{n}\:=\:\sum_{k=0}^n
\sbnm{u}{k}\sbnm{v}{n-k}.\]

\begin{exam}[$\lam=2$] For $\mu=0,1,2$, 
there hold partial fraction expansions
\xalignz{(n!)^{2-\mu}
\frac{(1-x)_n^{\mu}}{(x)^2_{n+1}}
=&\sum^n_{k=0}
\sbnm{n}{k}^2\sbnm{n+k}{k}^{\mu}
\bigg\{\frac{1}{(x+k)^2}\\
+&\frac{1}{x+k}
\Big[(2+\mu)H_k-2H_{n-k}-\mu H_{n+k}\Big]\bigg\} }
and the three corresponding harmonic number identities:
\[0\:=\:\sum^n_{k=0}
\sbnm{n}{k}^2\sbnm{n+k}{k}^{\mu}
\Big\{(2+\mu)H_k-2H_{n-k}-\mu H_{n+k}\Big\}.\]
\end{exam}

For $\mu=0$, the corresponding partial fraction
expansion has been given by \eqref{pf0/2}. 

\begin{exam}[$\lam=3$] For $\mu=0,1,2,3$, 
there hold partial fraction expansions
\xalignz{(n!)^{3-\mu}
\frac{(1-x)_n^{\mu}}{(x)^3_{n+1}}
=&\sum^n_{k=0}(-1)^k
\sbnm{n}{k}^3\sbnm{n+k}{k}^{\mu}
\bigg\{\frac{1}{(x+k)^3}\\
+&\frac{1}{(x+k)^2}
\Big[(3+\mu)H_k-3H_{n-k}-\mu H_{n+k}\Big]\\
+&\frac{1/2}{x+k}
\Bigg[\qdn\mult{rl}{\+\big\{(3+\mu)H_k\,-\,3H_{n-k}\,-\,\mu H_{n+k}\big\}^2\\[1mm]
+\+\big\{(3+\mu)H_k^{\ang{2}}+3H_{n-k}^{\ang{2}}-\mu H_{n+k}^{\ang{2}}\big\}}
\qdn\Bigg]\Bigg\} }
and the four corresponding harmonic number identities:
\[0\:=\:\sum^n_{k=0}(-1)^k
\sbnm{n}{k}^3\sbnm{n+k}{k}^{\mu}
\Bigg\{\qdn\mult{rl}{\+\big[(3+\mu)H_k\,-\,3H_{n-k}\,-\,\mu H_{n+k}\big]^2\\[1mm]
+\+\big[(3+\mu)H_k^{\ang{2}}+3H_{n-k}^{\ang{2}}-\mu H_{n+k}^{\ang{2}}\big]}
\qdn\Bigg\}.\]
\end{exam}

When $\mu=0$, the corresponding partial fraction decomposition 
and harmonic number identity have been exhibited respectively 
in (\ref{pf0/3a}-\ref{pf0/3b}) and \eqref{pf0/3c}.

\begin{exam}[$\lam=4$] For $\mu=0,1,2,3,4$, 
there hold partial fraction expansions
\xalignz{(n!)^{4-\mu}
\frac{(1\!-\!x)^\mu_n}{(x)^4_{n+1}}
&=\sum^n_{k=0}
\sbnm{n}{k}^4\sbnm{n+k}{k}^\mu
\Bigg\{\frac{1}{(x+k)^4}+
\frac{(4\!+\!\mu)H_k\!-4H_{n-k}\!-\mu H_{n+k}}{(x+k)^3}\\
&+\frac{1/2}{(x\!+\!k)^2}
\Big[{\sst\big\{(4+\mu)H_k-4H_{n-k}-\mu H_{n+k}\big\}^2
\:+\:\big\{(4+\mu)H_k^{\ang{2}}+4H^{\ang{2}}_{n-k}
	-\mu H^{\ang{2}}_{n+k}\big\}}\Big]\\
&+\frac{1/6}{x\!+\!k}\bigg[
\double{\big\{(4+\mu)H_k-4H_{n-k}-\mu H_{n+k}\big\}^3
+\,2\big\{(4+\mu)H_k^{\ang{3}}-4H^{\ang{3}}_{n-k}
	-\mu H^{\ang{3}}_{n+k}\big\}}
{+3\big\{(4+\mu)H_k-4H_{n-k}-\mu H_{n+k}\big\}
\times\big\{(4+\mu)H_k^{\ang{2}}+4H^{\ang{2}}_{n-k}
-\mu H^{\ang{2}}_{n+k}\big\}}\bigg]\Bigg\} }
and the five corresponding harmonic number identities:
\[0=\sum^n_{k=0}
\sbnm{n}{k}^4\sbnm{n+k}{k}^\mu
\bigg\{
\double{\big[(4+\mu)H_k-4H_{n-k}-\mu H_{n+k}\big]^3
+\,2\big[(4+\mu)H_k^{\ang{3}}-4H^{\ang{3}}_{n-k}
	-\mu H^{\ang{3}}_{n+k}\big]}
{+3\big[(4+\mu)H_k-4H_{n-k}-\mu H_{n+k}\big]
\times\big[(4+\mu)H_k^{\ang{2}}+4H^{\ang{2}}_{n-k}
-\mu H^{\ang{2}}_{n+k}\big]}\bigg\}.\]
\end{exam}

\begin{exam}[$\lam=5$] For $\mu=0,1,2,3,4,5$, 
there hold partial fraction expansions
\xalignz{(n!)^{5-\mu}
\frac{(1\!-\!x)^\mu_n}{(x)^5_{n+1}}
&=\sum^n_{k=0}(-1)^k
\sbnm{n}{k}^5\sbnm{n+k}{k}^\mu
\Bigg\{\frac{1}{(x+k)^5}
+\frac{\sst(5+\mu)H_k-5H_{n-k}-\mu H_{n+k}}{(x+k)^4}\\
&\qdn+\frac{1/2}{(x\!+\!k)^3}
\Big[{\sst\big\{(5+\mu)H_k-5H_{n-k}-\mu H_{n+k}\big\}^2
\:+\:\big\{(5+\mu)H_k^{\ang{2}}+5H^{\ang{2}}_{n-k}
	-\mu H^{\ang{2}}_{n+k}\big\}}\Big]\\
&\xqdn\!\qdp+\:\frac{1/6}{(x\!+\!k)^2}
\Bigg[\qdn\mult{r}{
\sst\big\{(5+\mu)H_k-5H_{n-k}-\mu H_{n+k}\big\}^3
+2\big\{(5+\mu)H_k^{\ang{3}}-5H^{\ang{3}}_{n-k}
	-\mu H^{\ang{3}}_{n+k}\big\}\\[1mm]
\sst+3\big\{(5+\mu)H_k-5H_{n-k}-\mu H_{n+k}\big\}
\times\big\{(5+\mu)H_k^{\ang{2}}+5H^{\ang{2}}_{n-k}
-\mu H^{\ang{2}}_{n+k}\big\}}\qdn\Bigg]\\
&\xqdn\qqdn+\left.
\frac{1/24}{x\!+\!k}\qdn
\left[\qdn\mult{c}{
{\sss\big\{(5+\mu)H_k-5H_{n-k}-\mu H_{n+k}\big\}^4
+3\big\{(5+\mu)H_k^{\ang{2}}+5H^{\ang{2}}_{n-k}
	-\mu H^{\ang{2}}_{n+k}\big\}^2
+6(5+\mu)H_k^{\ang{4}} }\\[1mm]
\sss+6\big\{(5+\mu)H_k-5H_{n-k}-\mu H_{n+k}\big\}^2
\times\big\{(5+\mu)H_k^{\ang{2}}+5H^{\ang{2}}_{n-k}
-\mu H^{\ang{2}}_{n+k}\big\}
+30H^{\ang{4}}_{n-k}\\[1mm]
\sss+8\big\{(5+\mu)H_k-5H_{n-k}-\mu H_{n+k}\big\}
\times\big\{(5+\mu)H_k^{\ang{3}}-5H^{\ang{3}}_{n-k}
	-\mu H^{\ang{3}}_{n+k}\big\}
-6\mu H^{\ang{4}}_{n+k}\;}\qdn\right]
\!\right\}  }
and the six corresponding harmonic number identities:
\xalignz{0\:&=\:\sum^n_{k=0}(-1)^k
\sbnm{n}{k}^5\sbnm{n+k}{k}^{\mu}\\
&\times\:
\left[\qdn\mult{c}{
{\sst\big\{(5+\mu)H_k-5H_{n-k}-\mu H_{n+k}\big\}^4
+3\big\{(5+\mu)H_k^{\ang{2}}+5H^{\ang{2}}_{n-k}
	-\mu H^{\ang{2}}_{n+k}\big\}^2
+6(5+\mu)H_k^{\ang{4}} }\\[1mm]
\sst+6\big\{(5+\mu)H_k-5H_{n-k}-\mu H_{n+k}\big\}^2
\times\big\{(5+\mu)H_k^{\ang{2}}+5H^{\ang{2}}_{n-k}
-\mu H^{\ang{2}}_{n+k}\big\}
+30H^{\ang{4}}_{n-k}\\[1mm]
\sst+8\big\{(5+\mu)H_k-5H_{n-k}-\mu H_{n+k}\big\}
\times\big\{(5+\mu)H_k^{\ang{3}}-5H^{\ang{3}}_{n-k}
	-\mu H^{\ang{3}}_{n+k}\big\}
-6\mu H^{\ang{4}}_{n+k}\;}\qdn\right].}
\end{exam}

By means of the standard partial fraction method, 
we can also derive the following algebraic identities
and the corresponding harmonic number formulae, even
though they are not consequences of the theorem 
and the corollary proved in the present paper.
\begin{exam} Partial fraction decomposition formula
\bnm
\frac{x(1-x)_n^2}{(x)^2_{n+1}}
&=&\frac{1}{x}+\sum^n_{k=1}
\sbnm{n}{k}^2\sbnm{n+k}{k}^2\\
&\times&\bigg\{\frac{-k}{(x+k)^2}
\:+\:\frac{1+2kH_{n+k}+2kH_{n-k}-4H_k}{x+k}\bigg\} 
\enm
and the corresponding harmonic number identity 
associated with Beukers' conjecture \emph{(cf.\:\cito{doron} 
and \cito{chu})}:
\[0\:=\:\sum^n_{k=1}
\sbnm{n}{k}^2\sbnm{n+k}{k}^2
\Big\{1+2kH_{n+k}+2kH_{n-k}-4H_k\Big\}.\]
\end{exam}
\begin{exam} Partial fraction decomposition formula
\bnm
\frac{(1-x)_n^2}{(1+x)^2_{n}}
&=&1+\sum^n_{k=1}
\sbnm{n}{k}^2\sbnm{n+k}{k}^2\\
&\times&\bigg\{\frac{k^2}{(x+k)^2}
\:-\:\frac{2k^2}{x+k}\Big(\frac{1}{k}+H_{n+k}+H_{n-k}
-2H_k\Big)\bigg\}
\enm
and the corresponding harmonic number identity: 
\[n(n+1)\:=\:\sum^n_{k=1}k^2\:
\sbnm{n}{k}^2\sbnm{n+k}{k}^2
\Big\{\frac{1}{k}+H_{n+k}+H_{n-k}-2H_k\Big\}.\]
\end{exam}

\begin{exam} Partial fraction expansion 
with denominator polynomial of type ``$2+1$''
\bnm
\frac{n!\times(2n)!}{(x)^2_{n+1}(1-x)_n}
\qdp&=&
\sum^n_{k=1}
\frac{\binm{2n}{n+k}}{\binm{n+k}{n+1}}
\frac{(-1)^{k}}{(1+n)(x-k)}\\
&+&\sum^n_{k=0}
\sbnm{n}{k}\sbnm{2n}{n+k}
\bigg\{\frac{1}{(x+k)^2}
+\frac{H_k+H_{n+k}-2H_{n-k}}
      {x+k}\bigg\}
\enm
and the corresponding harmonic number identity: 
\[\sum^n_{k=1}
\frac{(-1)^{k}}{1+n}
\frac{\binm{2n}{n+k}}{\binm{n+k}{n+1}}
\:=\:\sum^n_{k=0}
\sbnm{n}{k}\sbnm{2n}{n+k}
(2H_{n-k}-H_k-H_{n+k}).\]
\end{exam}

\begin{exam}  Partial fraction expansion 
with denominator polynomial of type ``$3+1$''
\bnm
\frac{(n!)^2\times(2n)!}{(x)^3_{n+1}(1-x)_n}
=\frac{1}{(1\!+\!n)^2}\sum^n_{k=1}
\frac{\binm{2n}{n+k}}{\binm{n+k}{n+1}^2}
\frac{(-1)^{k}}{x\!-\!k}
+\sum^n_{k=0}(-1)^k
\sbnm{n}{k}^2\sbnm{2n}{n+k}
\xquad\qqdn\zero\\
\times\:
\bigg\{\!\frac{1}{(x+k)^3}
+\frac{\sst2H_k+H_{n+k}-3H_{n-k}}
       {(x+k)^2}
+\frac{1/2}{x+k}\Big[\qdn\mult{c}
{\sst\big(2H_k\:+\:H_{n+k}\:-\:3H_{n-k}\big)^2\\
\sst+\big(2H_k^{\ang{2}}
+H_{n+k}^{\ang{2}}+3H_{n-k}^{\ang{2}}\big)}
\qdn\Big]\!\bigg\}\:\zero
\enm
and the corresponding harmonic number identity: 
\[\frac{-2}{(1+n)^2}
\sum^n_{k=1}(-1)^{k}
\frac{\binm{2n}{n+k}}{\binm{n+k}{n+1}^2}
=\sum^n_{k=0}(-1)^k
\sbnm{n}{k}^2\sbnm{2n}{n+k}
\bigg\{\qdn\mult{c}
{\sst\big(2H_k\:+\:H_{n+k}\:-\:3H_{n-k}\big)^2\\
\sst+\big(2H_k^{\ang{2}}
+H_{n+k}^{\ang{2}}
+3H_{n-k}^{\ang{2}}\big)}
\qdn\bigg\}.\]
\end{exam}

\begin{exam} Partial fraction decomposition formula
with denominator polynomial of type ``$3+2$''
\xalignz{
\frac{n!\times\big\{(2n)!\big\}^2}{(x)^3_{n+1}(1-x)^2_n}
&=\frac{1}{1+n}\sum^n_{k=1}
\frac{\binm{2n}{n+k}^2}{\binm{n+k}{n+1}}
\bigg\{\frac{1}{(x-k)^2}
+\frac{H_{k-1}+2H_{n-k}-3H_{n+k}}{x-k}\bigg\}\\
&+\sum^n_{k=0}(-1)^k
\sbnm{n}{k}\sbnm{2n}{n+k}^2
\bigg\{\frac{1}{(x+k)^3}
+\frac{\sst H_k+2H_{n+k}-3H_{n-k}}
      {(x+k)^2}\\
&+\frac{1}{2(x+k)}\Big[\sst
\big(H_k\:+\:2H_{n+k}\:-\:3H_{n-k}\big)^2
+\big(H_k^{\ang{2}}+2H_{n+k}^{\ang{2}}+3H_{n-k}^{\ang{2}}\big)
\Big]\bigg\}  }
and the corresponding harmonic number identity: 
\xalignz{
\sum^n_{k=0}&(-1)^k
\sbnm{n}{k}\sbnm{2n}{n+k}^2
\bigg\{\qdn\mult{c}
{\sst\big(H_k\:+\:2H_{n+k}\:-\:3H_{n-k}\big)^2\\
\sst+\big(H_k^{\ang{2}}+2H_{n+k}^{\ang{2}}+3H_{n-k}^{\ang{2}}\big)}
\qdn\bigg\}\\
=\qdp&\frac{2}{1+n}\sum^n_{k=1}
\frac{\binm{2n}{n+k}^2}{\binm{n+k}{n+1}}
\big(3H_{n+k}-H_{k-1}-2H_{n-k}\big).}
\end{exam}

\begin{exam} For nonnegative integer $\theta$ with 
$0\le\theta<4+4n$, there hold partial fraction 
decomposition formulae:
\bln\xalignx{
\frac{(n!)^4x^\theta}{(x)^4_{n+1}}
&=\sum^n_{k=0}
\sbnm{n}{k}^4
\bigg\{\frac{(-k)^{\theta}}{(x+k)^4}
+\frac{(-k)^{\theta-1}}{(x+k)^3}
\Big[\theta-4k(H_k-H_{n-k})\Big]\\
&+\frac{(-k)^{\theta-2}}{2(x+k)^2}
\Big[\sst\big\{\theta-4k(H_k-H_{n-k})\big\}^2
-\big\{\theta-4k^2(H_k^{(2)}+H_{n-k}^{(2)})\big\}\Big]\\
&+\frac{(-k)^{\theta-3}}{6(x+k)}
\bigg[\qdn\mult{r}{\sst
\big\{\theta\!-\!4k(H_k-H_{n-k})\big\}^3
+2\big\{\theta-4k^3(H_k^{(3)}\!-\!H_{n-k}^{(3)})\big\}\\
\sst-3\big\{\theta\!-\!4k(H_k-H_{n-k})\big\}
\big\{\theta\!-\!4k^2(H_k^{(2)}\!+\!H_{n-k}^{(2)})\big\}
\qdp}\qdn\bigg]\bigg\}.}\eln
\end{exam}
The corresponding harmonic identities read as 
\blm
&&\sum^n_{k=0}{k^{\theta-3}}
\sbnm{n}{k}^4
\bigg[\qdn\mult{r}{\sst
\big\{\theta\!-\!4k(H_k-H_{n-k})\big\}^3
+2\big\{\theta-4k^3(H_k^{(3)}\!-\!H_{n-k}^{(3)})\big\}\\
\sst-3\big\{\theta\!-\!4k(H_k-H_{n-k})\big\}
\big\{\theta\!-\!4k^2(H_k^{(2)}\!+\!H_{n-k}^{(2)})\big\}
\qdp}\qdn\bigg]\qquad\\
&&\;=\:\bigg\{\mult{cr}{
0,&0\le\theta\le2+4n\:\\
6(n!)^4,&\theta=3+4n.}
\elm 
For $\theta=0,1,2$, the corresponding results 
to this identity have been conjectured 
by Weideman~\citu{pade-1}{Eq\:21} and confirmed 
by Driver et al~\citu{pade-2}{Eq\:20}.
In particular, we recover, with the case $\theta=1$, 
the identity found by Driver et al~\citu{pade-3}{Eq\:21}.

The list can be endless. However, we are not bothered 
to extend it further. The interested reader can do
that for enjoyment.



\section*{Appendix A. Table for $\Ome$-Coefficients Computed 
		via \eqref{bell-k} and \eqref{leibniz-k}}
\renewcommand{\theequation}{A\arabic{equation}}
\setcounter{equation}{0}\vspace*{-9mm}
\xalignz{\Ome_0(\lam,\mu,-k)
&\equiv\tagg{1}{}1.\\
\Ome_1(\lam,\mu,-k)
&=\tagg{1}{}\lam\Big\{H_k-H_{n-k}\Big\}+\mu\Big\{H_k-H_{n+k}\Big\}.\\
\Ome_2(\lam,\mu,-k)
&=\tagg{1}{a}\Big\{\lam\big(H_k-H_{n-k}\big)
+\mu\big(H_k-H_{n+k}\big)\Big\}^2\\
&+\tagg{0}{b}{\lam}\Big\{H_k^{\ang{2}}+H_{n-k}^{\ang{2}}\Big\}
+{\mu}\Big\{H_k^{\ang{2}}-H_{n+k}^{\ang{2}}\Big\}.\quad}
\vspace*{-7mm}\xalignz{\Ome_3(\lam,\mu,-k)
&=\tagg{1}{a}\Big\{\lam\big(H_k-H_{n-k}\big)
+\mu\big(H_k-H_{n+k}\big)\Big\}^3\\
&+\tagg{0}{b}2
\Big\{\lam\big(H_k^{\ang{3}}-H_{n-k}^{\ang{3}}\big)
+\mu\big(H_k^{\ang{3}}-H_{n+k}^{\ang{3}}\big)\Big\}\\
&+\tagg{0}{c}3
\Big\{\lam\big(H_k-H_{n-k}\big)
+\mu\big(H_k-H_{n+k}\big)\Big\}\\
&\tagg{0}{d}\times
\Big\{\lam\big(H_k^{\ang{2}}+H_{n-k}^{\ang{2}}\big)
+\mu\big(H_k^{\ang{2}}-H_{n+k}^{\ang{2}}\big)\Big\}.}
\vspace*{-7mm}\xalignz{\Ome_4(\lam,\mu,-k)
&=\tagg{1}{a}
\Big\{\lam\big(H_k-H_{n-k}\big)
+\mu\big(H_k-H_{n+k}\big)\Big\}^4\\
&+\tagg{0}{b}6
\Big\{\lam\big(H_k^{\ang{4}}+H_{n-k}^{\ang{4}}\big)
+\mu\big(H_k^{\ang{4}}-H_{n+k}^{\ang{4}}\big)\Big\}\\
&+\tagg{0}{c}8
\Big\{\lam\big(H_k-H_{n-k}\big)
+\mu\big(H_k-H_{n+k}\big)\Big\}\\
&\tagg{0}{d}\times
\Big\{\lam\big(H_k^{\ang{3}}-H_{n-k}^{\ang{3}}\big)
+\mu\big(H_k^{\ang{3}}-H_{n+k}^{\ang{3}}\big)\Big\}\\
&+\tagg{0}{e}6
\Big\{\lam\big(H_k-H_{n-k}\big)
+\mu\big(H_k-H_{n+k}\big)\Big\}^2\\
&\tagg{0}{f}\times
\Big\{\lam\big(H_k^{\ang{2}}+H_{n-k}^{\ang{2}}\big)
+\mu\big(H_k^{\ang{2}}-H_{n+k}^{\ang{2}}\big)\Big\}\\
&+\tagg{0}{g}3
\Big\{\lam\big(H_k^{\ang{2}}+H_{n-k}^{\ang{2}}\big)
+\mu\big(H_k^{\ang{2}}-H_{n+k}^{\ang{2}}\big)\Big\}^2.}
\vspace*{-7mm}\xalignz{\Ome_5(\lam,\mu,-k)
&=\tagg{1}{a}
\Big\{\lam\big(H_k-H_{n-k}\big)
+\mu\big(H_k-H_{n+k}\big)\Big\}^5\\
&+\tagg{0}{b}24
\Big\{\lam\big(H_k^{\ang{5}}-H_{n-k}^{\ang{5}}\big)
+\mu\big(H_k^{\ang{5}}-H_{n+k}^{\ang{5}}\big)\Big\}\\
&+\tagg{0}{c}10
\Big\{\lam\big(H_k-H_{n-k}\big)
+\mu\big(H_k-H_{n+k}\big)\Big\}^3\\
&\tagg{0}{d}\times
\Big\{\lam\big(H_k^{\ang{2}}+H_{n-k}^{\ang{2}}\big)
+\mu\big(H_k^{\ang{2}}-H_{n+k}^{\ang{2}}\big)\Big\}\\
&+\tagg{0}{e}20
\Big\{\lam\big(H_k-H_{n-k}\big)
+\mu\big(H_k-H_{n+k}\big)\Big\}^2\\
&\tagg{0}{f}\times
\Big\{\lam\big(H_k^{\ang{3}}-H_{n-k}^{\ang{3}}\big)
+\mu\big(H_k^{\ang{3}}-H_{n+k}^{\ang{3}}\big)\Big\}\\
&+\tagg{0}{g}15
\Big\{\lam\big(H_k-H_{n-k}\big)
+\mu\big(H_k-H_{n+k}\big)\Big\}\\
&\tagg{0}{h}\times
\Big\{\lam\big(H_k^{\ang{2}}+H_{n-k}^{\ang{2}}\big)
+\mu\big(H_k^{\ang{2}}-H_{n+k}^{\ang{2}}\big)\Big\}^2\\
&+\tagg{0}{i}30
\Big\{\lam\big(H_k-H_{n-k}\big)
+\mu\big(H_k-H_{n+k}\big)\Big\}\\
&\tagg{0}{j}\times
\Big\{\lam\big(H_k^{\ang{4}}+H_{n-k}^{\ang{4}}\big)
+\mu\big(H_k^{\ang{4}}-H_{n+k}^{\ang{4}}\big)\Big\}\\
&+\tagg{0}{k}20
\Big\{\lam\big(H_k^{\ang{2}}+H_{n-k}^{\ang{2}}\big)
+\mu\big(H_k^{\ang{2}}-H_{n+k}^{\ang{2}}\big)\Big\}\\
&\tagg{0}{l}\times
\Big\{\lam\big(H_k^{\ang{3}}-H_{n-k}^{\ang{3}}\big)
+\mu\big(H_k^{\ang{3}}-H_{n+k}^{\ang{3}}\big)\Big\}.}


\section*{Appendix B. Table for $\varpi$-Coefficients 
		Computed via \eqref{bell-ka}}
\renewcommand{\theequation}{B\arabic{equation}}
\setcounter{equation}{0}\vspace*{-9mm}
\xalignz{\varpi_0(\lam,-k)
&\equiv\tagg{1}{}1.\\
\varpi_1(\lam,-k)
&=\tagg{1}{}
\lam\Big\{H_k-H_{n-k}\Big\}.\\
\varpi_2(\lam,-k)
&=\tagg{1}{}{\lam^2}
\Big\{H_k-H_{n-k}\Big\}^2
+{\lam}
\Big\{H_k^{(2)}+H_{n-k}^{(2)}\Big\}.\qdp}
\vspace*{-7mm}\xalignz{\varpi_3(\lam,-k)
&=\tagg{1}{a}{\lam^3}
\Big\{H_k-H_{n-k}\Big\}^3
+2{\lam}
\Big\{H_k^{(3)}-H_{n-k}^{(3)}\Big\}\\
&+\tagg{0}{b}3{\lam^2}
\Big\{H_k-H_{n-k}\Big\}\times
\Big\{H_k^{(2)}+H_{n-k}^{(2)}\Big\}.}
\vspace*{-7mm}\xalignz{\varpi_4(\lam,-k)
&=\tagg{1}{a}{\lam^4}
\Big\{H_k-H_{n-k}\Big\}^4
+6{\lam}
\Big\{H_k^{(4)}+H_{n-k}^{(4)}\Big\}\\
&+\tagg{0}{b}8{\lam^2}
\Big\{H_k-H_{n-k}\Big\}\times
\Big\{H_k^{(3)}-H_{n-k}^{(3)}\Big\}\\
&+\tagg{0}{c}6{\lam^3}
\Big\{H_k-H_{n-k}\Big\}^2\times
\Big\{H_k^{(2)}+H_{n-k}^{(2)}\Big\}\\
&+\tagg{0}{d}3{\lam^2}
\Big\{H_k^{(2)}+H_{n-k}^{(2)}\Big\}^2.}
\vspace*{-7mm}\xalignz{\varpi_5(\lam,-k)
&=\tagg{1}{a}{\lam^5}
\Big\{H_k-H_{n-k}\Big\}^5
+24{\lam}
\Big\{H_k^{(5)}-H_{n-k}^{(5)}\Big\}\\
&+\tagg{0}{b}10{\lam^4}
\Big\{H_k-H_{n-k}\Big\}^3\times
\Big\{H_k^{(2)}+H_{n-k}^{(2)}\Big\}\\
&+\tagg{0}{c}20{\lam^3}
\Big\{H_k-H_{n-k}\Big\}^2\times
\Big\{H_k^{(3)}-H_{n-k}^{(3)}\Big\}\\
&+\tagg{0}{d}15{\lam^3}
\Big\{H_k-H_{n-k}\Big\}\times
\Big\{H_k^{(2)}+H_{n-k}^{(2)}\Big\}^2\\
&+\tagg{0}{e}30{\lam^2}
\Big\{H_k-H_{n-k}\Big\}\times
\Big\{H_k^{(4)}+H_{n-k}^{(4)}\Big\}\\
&+\tagg{0}{f}20{\lam^2}
\Big\{H_k^{(2)}+H_{n-k}^{(2)}\Big\}\times
\Big\{H_k^{(3)}-H_{n-k}^{(3)}\Big\}.}


\section*{Appendix C. Table for $\ome$-Coefficients 
		Computed via \eqref{bell-kb}}
\renewcommand{\theequation}{C\arabic{equation}}
\setcounter{equation}{0}\vspace*{-9mm}
\xalignz{\ome_0(\mu,-k)
&\equiv\tagg{1}{}1.\\
\ome_1(\mu,-k)
&=\tagg{1}{}\mu\Big\{H_k-H_{n+k}\Big\}.\\
\ome_2(\mu,-k)
&=\tagg{1}{}{\mu^2}
\Big\{H_k-H_{n+k}\Big\}^2
+{\mu}
\Big\{H_k^{\ang{2}}-H_{n+k}^{\ang{2}}\Big\}.}
\vspace*{-7mm}\xalignz{\ome_3(\mu,-k)
&=\tagg{1}{a}{\mu^3}
\Big\{H_k-H_{n+k}\Big\}^3
+2{\mu}
\Big\{H_k^{\ang{3}}-H_{n+k}^{\ang{3}}\Big\}\\
&+\tagg{0}{b}3{\mu^2}
\Big\{H_k-H_{n+k}\Big\}\times
\Big\{H_k^{\ang{2}}-H_{n+k}^{\ang{2}}\Big\}.}
\vspace*{-7mm}\xalignz{\ome_4(\mu,-k)
&=\tagg{1}{a}{\mu^4}
\Big\{H_k-H_{n+k}\Big\}^4
+6{\mu}
\Big\{H_k^{\ang{4}}-H_{n+k}^{\ang{4}}\Big\}\\
&+\tagg{0}{b}8{\mu^2}
\Big\{H_k-H_{n+k}\Big\}\times
\Big\{H_k^{\ang{3}}-H_{n+k}^{\ang{3}}\Big\}\\
&+\tagg{0}{c}6{\mu^3}
\Big\{H_k-H_{n+k}\Big\}^2\times
\Big\{H_k^{\ang{2}}-H_{n+k}^{\ang{2}}\Big\}\\
&+\tagg{0}{d}3{\mu^2}
\Big\{H_k^{\ang{2}}-H_{n+k}^{\ang{2}}\Big\}^2.}
\vspace*{-7mm}\xalignz{\ome_5(\mu,-k)
&=\tagg{1}{a}{\mu^5}
\Big\{H_k-H_{n+k}\Big\}^5
+24{\mu}
\Big\{H_k^{\ang{5}}-H_{n+k}^{\ang{5}}\Big\}\\
&+\tagg{0}{b}10{\mu^4}
\Big\{H_k-H_{n+k}\Big\}^3\times
\Big\{H_k^{\ang{2}}-H_{n+k}^{\ang{2}}\Big\}\\
&+\tagg{0}{c}20{\mu^3}
\Big\{H_k-H_{n+k}\Big\}^2\times
\Big\{H_k^{\ang{3}}-H_{n+k}^{\ang{3}}\Big\}\\
&+\tagg{0}{d}15{\mu^3}
\Big\{H_k-H_{n+k}\Big\}\times
\Big\{H_k^{\ang{2}}-H_{n+k}^{\ang{2}}\Big\}^2\\
&+\tagg{0}{e}30{\mu^2}
\Big\{H_k-H_{n+k}\Big\}\times
\Big\{H_k^{\ang{4}}-H_{n+k}^{\ang{4}}\Big\}\\
&+\tagg{0}{f}20{\mu^2}
\Big\{H_k^{\ang{2}}-H_{n+k}^{\ang{2}}\Big\}\times
\Big\{H_k^{\ang{3}}-H_{n+k}^{\ang{3}}\Big\}.}

\end{document}